\documentclass[a4paper,12pt]{article}

\usepackage{amsmath, latexsym, amsfonts, amssymb, amsthm, amscd}

\newtheorem{theorem}{Theorem}

\newtheorem{lemma}{Lemma}
\newtheorem{rem}{Remark}

\newtheorem{exple}{Example}

\newcommand{\R}{\mbox{\rm I\hspace{-0.02in}R}}
\newcommand{\F}{\mbox{\rm I\hspace{-0.02in}F}}

\def\QED{\hfill\vrule height 1.5ex width 1.4ex depth -.1ex \vskip20pt}

\begin{document}

\title{On continuity properties for option prices in exponential L\'evy models}
\maketitle

\begin{center}
{\large S. Cawston}\footnote{$^{,2}$ LAREMA, D\'epartement de
Math\'ematiques, Universit\'e d'Angers, 2, Bd Lavoisier - 49045,
\\\hspace*{.4in}{\sc Angers Cedex 01.}

\hspace*{.05in}$^1$E-mail: suzanne.cawston@univ-angers.fr$\;\;\;$
$^2$E-mail: lioudmila.vostrikova@univ-angers.fr}{\large
 and  L. Vostrikova$^2$}
\end{center}
\vspace{0.2in}

\begin{abstract} 
For a converging sequence of exponential L\'evy models, we give
conditions under which the associated sequence of option prices
converges. We also study the behaviour of the prices when no such
convergence holds. We then consider two special cases, first when the martingale
measure is chosen by minimisation of entropy and then when it minimises Hellinger integrals.


\end{abstract}

\section{Introduction and main results}
One of the most famous models for the price of a risky asset is geometric Brownian motion.   Although its simplicity is appealing, the Black-Scholes model isn't very
accurate for a number of reasons. One of these is that the fit for the
law of logreturns of stock-prices is relatively poor (cf. \cite{EK}) : normal distributions
don't give enough mass either around zero or to the tails. Therefore, a number
of other models have been developed which take these features into account. 
\par Since economic time doesn't follow the natural time scale but
a 'financial clock' which can be represented by a random process
$(\tau_t)_{t\geq 0}$, the basic idea consists in using time
changes. This is for example the case for hyperbolic models which
were introduced by O.E. Barndorff-Nielsen in the context of modelling
the size of sand deposits, before being used in finance. In Generalized
Hyperbolic models GH$(\lambda,\alpha,\beta,\delta,\mu)$(see
\cite{EK},\cite{E},\cite{P}) the logarithm
of the stock price is assumed to be given by a process
$$X_t=\mu t+\beta\tau_t+W_{\tau_t}$$
where $\mu$ is a deterministic drift, $W$ is a Brownian motion and $\tau$ is a Generalized Inverse Gaussian process GIG$(\lambda,\delta,\sqrt{\alpha^2-\beta^2})$(see
\cite{BNH}) which is independent of 
$W$.  Just as in the Black-Scholes case, $X$ is a Levy process, i.e
it has stationary and independent increments. However, the paths of $X$
are no longer continuous. In fact, it can be shown that $X$ has no
continuous martingale component, but has an infinite number of jumps and infinite
variation over any time-interval.
\par The main interest of this vast
family of models is that it allows for an excellent fit both for daily
log-returns, and for intraday data (cf. \cite{EF}). Furthermore, is
was shown in \cite{EH} that a number of
other popular models, including the Black-Scholes case, can be obtained as
limiting cases. In this paper, we will focus in our examples on normal inverse
gaussian (NIG) processes which correspond to the case
$\lambda=-1/2$. These processes preserve the goodness of fit
properties of more general GH models, but have the extra advantage that the form of the law of $X_t$
is the same at all times $t$ which helps speed up simulations and numerical pricing. The Levy measure of NIG$(\alpha,\beta,\delta,\mu)$ can also be expressed
relatively simply under the form
\begin{equation}\label{levnig}
\nu(dx)=\frac{\alpha\delta}{\pi}\frac{e^{\beta
    x}}{|x|}K_1(\alpha|x|)dx
\end{equation}
where $\alpha , \beta , \delta$ are parameters, and $K_1$ denotes the modified Bessel function of the third kind of parameter 1.    
\par Despite their qualities, one drawback of GH 
models is that $X$ necessarily has infinite variation whereas it was
suggested in \cite{CGMY2} that this isn't always appropriate for modelling
financial data. One family of models which allows for both
finite and infinite variation is obtained when assuming $X$ is a CGMY process
(cf. \cite{CGMY2}), that is a Levy process with no continuous martingale
component and with a Levy measure $\nu$ of the form 
\begin{equation}\label{levcgmy}
\nu(dx)=\frac{C}{|x|^{Y+1}}(e^{-G|x|}I_{\{x<0\}}+e^{-M|x|}I_{\{x>0\}})dx
\end{equation}
where $ C, G, M, Y$ are the parameters of the model.
In fact, this family generalises the Variance-Gamma process introduced
in \cite{MS},
which corresponds to the case $Y=1$ and can be obtained in the same
way as GH
processes by subordinating a Brownian Motion, but this time using a Gamma clock.
  
\par All these exponential Levy models depend on several parameters. For example, in the
Black-Scholes model, the parameters are the drift $\mu$ and the
volatility $\sigma$, we have four parameters in a CGMY model and five for a GH model.
 In practice, these are
usually calibrated and assumed to be constant over some interval of
time. However, as the information which is available increases continuously, it is important
to consider a dynamic approach, in which we have a sequence of parameters  and the corresponding sequence of stochastic processes. As far
as we know, this kind of problem was first considered by A.N. Shiryaev in the
case of the convergence of a sequence of Cox-Ross-Rubinstein models to
a Black-Scholes model (cf. \cite{SH}, Chapter 6.3d).
\par In its simplest form, this approach leads to the following formalisation. Assume we are given $(\Omega^n,\mathcal{F}^n,\mathbb{F}^n,P^n)_{n\geq 1}$, a sequence of stochastic bases with a right-continuous filtration $\mathbb{F}^n=(\mathcal{F}_t^n)_{t\geq 0}$ which is completed with respect to $P^n$ and with $\mathcal{F}_n=\bigvee_{t\geq 0}\mathcal{F}_t^n$.  
\par Assume that $(X^n)_{n\geq 1}$ is a sequence of Levy processes
with characteristics $(b_n,c_n,\nu_n)$ respectively, (see
\cite{B},\cite{Sa}).  We recall that $b_n$ represents the drift, $c_n$ the quadratic
variation and $\nu_n$ the Levy measure which satisfies the usual
condition
\begin{equation}
\label{1}
\int _{\mathbb{R}^*}(x^2\wedge 1)\nu _n(dx) < \infty .
\end{equation}
The characteristic function of $X^n_t$ for $t\geq 0$ and $u\in\mathbb{R}$ is equal to 
$$\phi^n_t(u)=E_{P^n}exp(iuX_t^n)=exp(t\psi_n(iu))$$
where $E_{P^n}$ is the expectation with respect to $P^n$ and $\psi_n(u)$ is the characteristic exponent given by 
$$\psi_n(u)=b_nu+\frac{c_n}{2}u^2+\int_{\mathbb{R}^*}[e^{ux}-1-uh(x)]\nu_n(dx)$$
In the last formula, $h$ denotes the truncation function.

We assume that we are given a risky asset which is modelled by 
\begin{equation}
\label{2}
S^n_t = S^n_0\exp\{X ^n_t\}
\end{equation}
and a non-risky asset 
\begin{equation}
\label{3}
B^n_t = B^n_0\exp\{r _n t\},
\end{equation}
where $r_n$ is the interest rate, $r_n\geq 0$. We also consider an option of maturity $T>0$. Let $g$ denote the associated payoff function which we assume to be continuous on the Skorokhod space $D[0,T]$ and to satisfy
\begin{equation}
\label{4}
g(Y)\leq  A \sup_{0\leq t \leq T} | Y_t| +B
\end{equation}
for $Y\in D[0,T]$, where $A,B$ are positive constants. For european call options, we have $g(Y)=(Y_{T}-K)^+$, for asian options $g(Y)=(Y_{T}-\frac{1}{T}\int_{0}^T Y_t dt)^+$ and for lookback options $g(Y)=(Y_{T}-\alpha \inf_{0\leq t\leq T}Y_{t})^+$, with $\alpha>1$. All these options, as well as the corresponding put options, satisfy (\ref{4}). 
\par We assume that for every $n\geq 1$ the set of equivalent
martingale measures is not empty and that we have chosen for some
reason one of these measures, say $Q^n$. The selection of the
martingale measure when it is not unique is an important step in the
derivation of an option price. Several approaches have been developed
in the literature, for instance minimisation of entropy \cite{MIF},
\cite{MI},\cite{ES}, of a Hellinger distance \cite{CS},\cite{CSL} or
of a Hellinger integral of order $q$ \cite{JKM}. From an economical
point of view, such approaches are motivated by their link with the
dual problem of utility maximisation
\cite{KrS},\cite{BF},\cite{F},\cite{GR},\cite{K} or mean-squared risk
minimisation \cite{FS}, \cite{S1},\cite{S2}.
\par Assume the equivalent martingale measure has been chosen and is equal to $Q^n$. Then the option price is equal to 
$$\mathbb{C}_{T}^n=E_{Q^n}[g(\tilde{S}^n)]$$
where $\tilde{S}^n_t=S_t^n/B_t^n$ is the discounted price of the risky asset. 
\par We assume that as information increases, our model tends in some sense to a limiting model which is related to a Levy process $X=(X_{t})_{t\geq 0}$ given on a canonical basis $(\Omega,\mathcal{F},\mathbb{F},P)$ and with characteristics $(b,c,\nu)$. This model consists again of two assets 
\begin{equation}
\label{5a}
S_t= S_0 \exp\{X_t\}
\end{equation}
\begin{equation}
\label{6}
B_t= B_0 \exp\{r t\}
\end{equation}
where $r\geq 0$ is the interest rate. This approach leads to several natural questions. First of all, under what conditions does an equivalent martingale measure exist for the limiting model ? If we assume for simplicity that such a measure exists and is equal to $Q$, we set 
$$\mathbb C _T= E_Q[g(\tilde{S})]$$
where $\tilde{S}_t=S_t/B_t$ is as before the corresponding discounted price. The next question is what conditions then ensure convergence of the option prices : 
\begin{equation}
\label{5}
\lim _{n\rightarrow \infty} \mathbb C^n_T = \mathbb C _T.
\end{equation}

In this article, we give sufficient conditions which ensure the
existence of $\lim_{n\to+\infty}\mathbb{C}_{T}^n$ as well as
conditions for the existence of a martingale measure in the limiting
model and conditions for (\ref{5}). We also give examples of cases
when (\ref{5}) is not satisfied.
\par Let $(\beta^n,Y^n)$ and $(\beta,Y)$ be the parameters which arise
in the Girsanov theorem and determine the density processes of the
changes of measures from $P^n$ and $P$ to the equivalent martingale
measures $Q^n$ and $Q$. From now on, we will refer to these parameters
as the Girsanov parameters of a change of measure. We assume that
under the new martingale measures, the processes $X^n$ and $X$ remain
Levy processes. This assumption isn't too restrictive for the pricing
of options as it has been shown in \cite{EJ},\cite{J} that for a
number of standard models the price interval covered by these
structure preserving measures is the whole non-arbitrage
interval. Since the measures $Q^n$ and $Q$ are martingale measures, we
must have (cf. \cite{JS} p. 556)
$$\int_{x>1}(e^x-1)Y^n(x)\nu_n(dx)<+\infty \text{ , }\int_{x>1}(e^x-1)Y(x)\nu(dx)<+\infty$$
\begin{theorem}\label{t1}
Assume the payoff function $g$ satisfies condition (\ref{4}). Assume furthermore that
\begin{enumerate}
\item $\lim _{n\rightarrow \infty} S^n_0 = S_0,$
\item $\lim _{n\rightarrow \infty} \left(  c_n + \int_{\mathbb{R}^*} h^2(x) Y^n(x) \nu_n(dx) \right) =  c + \int_{\mathbb{R}^*} h^2(x) Y(x) \nu (dx) $
\item $\lim _{n\rightarrow \infty} \int _{\mathbb{R}^*} (e^x -1)f(x) Y^n(x) \nu_n(dx) = \int _{\mathbb{R}^*} (e^x -1)  f(x) Y(x) \nu(dx)$
\end{enumerate}
for all continuous bounded functions $f$ which satisfy the condition $\lim_{x\to 0}\frac{f(x)}{x}=0$. 
\par Then we have convergence (\ref{5}) for option prices. 
\end{theorem}
\begin{rem} The main point in the proof of Theorem \ref{t1} consists in establishing the uniform integrability of the family $(\sup_{0\leq t\leq T}S^n_t)_{n\geq 1}$. This is shown in Lemma \ref{L3} of section 2 and is based on the Wiener-Hopf factorisation. 
\end{rem}

Next, we assume that we have additional information, namely that 
\begin{equation}\label{a}
\mathcal{L}(S^n|P^n)\longrightarrow \mathcal{L}(S|P)
\end{equation}
The question is then to know what conditions need to be added to (\ref{a}) to ensure that
\begin{equation}\label{c}
\mathcal{L}(S^n|Q^n)\longrightarrow \mathcal{L}(S|P^*)
\end{equation}
where $Q^n$ is a martingale measure for $S^n$ and $P^*$ is some measure which is absolutely continuous with respect to $P$. In general, the answer is known (cf \cite{JS}) and linked to the convergence 
\begin{equation}\label{d}
\mathcal{L}((Z^n,S^n)|P^n)\longrightarrow \mathcal{L}((Z,S)|P)
\end{equation}
where $Z^n$ and $Z$ denote the Radon-Nikodym derivatives of $Q^n$ and
$P^*$ with respect to $P^n$ and $P$ respectively. But in the case of
Levy processes, it will be easier to show (\ref{c}) by directly considering the convergence of the characteristics of the processes involved. 
\par We assume that the Girsanov parameters $(\beta_n,Y_n)$ and $(\beta,Y)$ associated with the changes of measures $P^n$ to $Q^n$ and $P$ to $P^*$ respectively, satisfy the following property in a neighbourhood of $0$ : 
\begin{equation}\label{e}
Y^n(x)=1+\beta^n x+o(x),
\end{equation}
\begin{equation}\label{f}
Y(x)=1+\beta x+o(x)
\end{equation}
where $o(x)$ in (\ref{e}) is uniform in $n$. 

\begin{theorem}\label{t1a}
 We assume that (\ref{a}),  (\ref{e}) and (\ref{f}) are satisfied. We assume furthermore that 
\par i) $\lim_{n\to+\infty}\beta^n=\beta$,
\par ii) For all $\epsilon>0$, 
$$\lim_{n\to+\infty}\int_{|x|\geq \epsilon}Y^n(x)\nu_n(dx)=\int_{|x|\geq \epsilon}Y(x)\nu(dx)$$
Then for all continuous bounded payoff function $g$, we have 
$$\lim_{n\to+\infty}\mathbb{C}_{T}^n=E_{P^*}[g(S)]$$
where $P^*$ is absolutely continuous with respect to $P$ and given by the Girsanov parameters $(\beta,Y)$. Moreover, $P^*$ is an equivalent martingale measure for $S$ if and only if the following conditions are satisfied : $Y>0 (\nu$-a.s.) and 
$$
\lim_{n\to+\infty}\int_{x\geq 1}e^xY^n(x)\nu_n(dx)=\int_{x\geq1}e^xY(x)\nu(dx)$$
Under this extra assumption, the prices of options converge for all payoff functions $g$ which satisfy (\ref{4}). 
\end{theorem}

In the sequel, we consider two particular choices of equivalent martingale measures, which minimise relative entropy and Hellinger integrals respectively. 
\par Let $Q$ and $P$ be two equivalent martingale measures. Then the entropy of $Q$ with respect to $P$ (or Kullback-Leibler information of $Q$ with respect to $P$) is given by 
$$H(Q|P)=E_{Q}(\ln (\frac{dQ}{dP}))=E_{P}(\frac{dQ}{dP}\ln (\frac{dQ}{dP}))$$
We recall that the martingale measure with minimal entropy is a
measure $P^{ME}$ such that the process $(exp(-rt)S_{t})_{t\geq 0}$ is
a martingale under $P^{ME}$ and such that for every equivalent martingale measure $Q$ 
$$H(P^{ME}|P)\leq H(Q|P)$$
It turns out (cf. \cite{KS},\cite{HS}) that in the case of Levy processes, if $P^{ME}$ exists, it is nothing else than an Esscher measure for the process $(\hat{X}_{t})_{t\geq 0}$ such that 
$$S_{t}=S_{0}\mathcal{E}(\hat{X})_{t}$$
where $\mathcal{E}(.)$ denotes the Doleans-Dade exponential. It is well know that $\hat{X}$ is a Levy process whose characteristic exponent is given by the formula (cf. \cite{MIF})
$$\hat{\psi}(u)=(b+\frac{1}{2}c)u+c\frac{u^2}{2}+\int_{\mathbb{R}^*}[e^{u(e^x-1)}-1-uh(x)]\nu(dx)$$
Let 
$$D=\{u\in\mathbb{R} : E_{P}(exp(u\hat{X}_1))<+\infty\}$$
We introduce the Esscher measure $P^u$ corresponding to $\hat{X}$ and $u\in D$ : 
$$\frac{d P_t^{u}}{d P_t}= \frac{\exp(u\hat{X}_t)}{E_P(\exp (u\hat{X}_t))}$$
where $P^u_t$ and $P_t$ denote the restriction of $P^u$ and $P$ to the $\sigma-$algebra $\mathcal{F}_{t}$, $t\geq 0$. In order for $P^u$ to be an equivalent martingale measure, the parameter $u$ should be equal to $\theta$, where $\theta$ satisfies the equation 
\begin{equation}
\label{8}
b+(\frac{1}{2}+\theta) c+ \int_{\mathbb{R}^*}\left[(e ^x -1)e^{\theta (e^x-1)} -h(x)\right]\nu(dx) = r
\end{equation}
The same results can be written for the measure $Q^n=(P^n)^{ME}$. Namely, the corresponding parameter of the Esscher measure satisfies the equation 
\begin{equation}
\label{9}
b_n+(\frac{1}{2}+\theta) c_n+ \int_{\mathbb{R}^*}[(e ^x -1)e^{\theta (e^x-1)} -h(x)]\nu_n(dx)= r_n
\end{equation}
We will assume that for every $n\geq 1$, equation (\ref{9}) has a
solution. This implies in particular that we exclude the case of
monotone Levy processes (cf Lemma \ref{L31} section 3.). In \cite{KS}, it was shown that under these
assumptions, the solution to (\ref{9}) is unique. But as we will see
later, this does not necessarily imply the existence of a solution to equation (\ref{8}). 
\par We define the set  
\begin{equation}
\label{u}
U= \{u\in \mathbb{R} : \overline{\lim}_{n\rightarrow \infty} \int _{x>1} e^{u(e^x -1)} \nu_n(dx) < \infty \}
\end{equation}
$U$ is the set on which the integrals considered are uniformly bounded for large values of $n$. It may be an open interval $]-\infty,\alpha[$ or a closed interval $]-\infty,\alpha]$ where $\alpha=\sup \{u:u\in U\}$. 

\begin{theorem}
\label{t2}
Assume the payoff function $g$ satisfies (\ref{4}). Assume furthermore that :
\begin{enumerate}
\item $\lim _{n\rightarrow \infty} S^n_0 = S_0,$
\item $\lim _{n\rightarrow \infty}b_n -r_n = b - r,$
\item $\lim _{n\rightarrow \infty} \left(  c_n + \int_{\mathbb{R}^*} h^2(x) \nu_n(dx) \right) =  c + \int_{\mathbb{R}^*} h^2(x)  \nu (dx) $
\item $\lim _{n\rightarrow \infty} \int _{\mathbb{R}^*} f(x)  \nu_n(dx) = \int _{\mathbb{R}^*}  f(x)  \nu(dx)$
\end{enumerate}
for all continuous bounded functions $f$ such that $\lim_{x\to 0}\frac{f(x)}{x^2}=0$. 
\par Then, if $\lim_{u\to\alpha^-}\hat{\psi}'(u)\geq 0$, there exists an equivalent martingale measure with minimal entropy for the limiting model. Furthermore, if $\lim_{u\to\alpha^-}\hat{\psi}'(u)>0$, we have convergence (\ref{5}) for the option prices and if $\lim_{u\to\alpha^-}\hat{\psi}'(u)=0$, (\ref{5}) holds at least for a subsequence. 
\par If $\lim_{u\to\alpha^-}\hat{\psi}'(u)<0$, the existence of a minimal equivalent martingale measure is not guaranteed, but if $g$ is bounded, we have 
$$\lim_{n'\to+\infty}\mathbb{C}_{T}^{n'}=E_{P^*} [g(S)]$$
where $n'$ is some subsequence and $P^*$ is a measure equivalent to $P$ which is not a martingale measure and whose Girsanov parameters are $(\alpha,exp(\alpha(e^x-1)))$. 
\end{theorem}
\begin{rem}
In the case when $\lim_{u\to\alpha^-}\hat{\psi}'(u)<0$, a martingale measure of minimal entropy may (Example 1, section 3.) or may not exist (Example 2, section 3.). Even if it does, the limit for the option prices is in general not $\mathbb{C}_{T}$. Moreover, using the Wiener-Hopf factorisation, we can show that 
$$\lim_{n\to+\infty}E_{Q}^n[\sup_{0\leq t\leq T}S_{t}^n]\neq E_{P^*}[\sup_{0\leq t\leq T}S_{t}]$$
since $E_{P^*}S_{T}\neq 1$. This means that the family of random
variables $(\sup_{0\leq t\leq T}S^n_{t})_{n\geq 1}$ is not uniformly
integrable and neither is $(S^n_T)_{n\geq 1}$. Consequently, we cannot extend the result to unbounded payoff functions. However, in some cases we can not only prove convergence but also derive an explicit expression for the limit. For instance, for european call options, we have 
$$\lim_{n\rightarrow \infty} \mathbb C^n_T = \lim_{n\rightarrow \infty} E_{Q^n}(S^n_T - K)^+ = E_{P^*} (S_T -K)^+ + 1- E_{P^*} (S_T).$$
\end{rem}\rm

An alternative choice for the martingale measure is related to the
$f^q$-martingales introduced in \cite{JKM}. These measures are special
cases of measures which minimise an $f$-divergence of Ciszar
\cite{CI}. Let $f$ be a convex function defined on $\mathbb{R}^{+,*}$
and $Q<<P$ be two probability measures on ($\Omega,\mathcal{F}$). The
$f$-divergence of $Q$ with respect to $P$ is then defined by 
$$f(Q|P)=E_{P}[f(\frac{dQ}{dP})]$$ 
It is easy to see that if $f(x)=x\ln x$ we obtain the Kullback-Leibler information or entropy of $Q$ with respect to $P$. If $f(x)=|x-1|$, we obtain the variational distance between $Q$ and $P$, for $f(x)=(1-x)^2$, we have quadratic variation and if $f(x)=(1-\sqrt{x})^2$, we have a Hellinger-distance, and finally if 
$$f(x)= \left\{
\begin{array}{ccc}
-x^q,&\mbox{ if } & 0<q<1,\\
x ^q, & \mbox{ if } & q>1 \mbox{ or } q<0,
\end{array}\right.$$
we obtain, up to a sign, the Hellinger integrals corresponding to $Q$ and $P$. The equivalent martingale measure which minimises these integrals has been studied in \cite{JKM} for $q>1$ and $q<0$, and in \cite{CSL} for $0<q<1$.

\par If we exclude monotone Levy processes and consider the larger set
of martingale measures which are only absolutely continuous with
respect to $P^n$, the $f^q$-martingale measure $Q^n$ exists and its
Girsanov parameters $(\beta_n,Y_n)$ are the solution to a minimisation
problem (see Lemma \ref{min}). $Q^n$ will then be equivalent to $P^n$
if and only if $Y^n>0$ $\nu_n$-a.s and we will assume in the sequel
that this condition holds. In order to exclude a trivial case, we will
also assume that $P^n$ is not a martingale measure itself. 
\par In the sequel, we will assume that the limiting process $X$ is not a monotone Levy process and that $\nu(\{m\})=\nu(\{M\})=0$ where $m$ and $M$ are the infimum and supremum of $supp(\nu)$.  
We introduce the integrals 
+$$I_n(q)=\int_{x\geq 1}e^{\frac{qx}{q-1}}\nu_n(dx)$$
which we assume to be finite. We define in the same way $I(q)$ where $\nu_n$ is replaced by $\nu$, but $I(q)$ may or may not be finite.

\begin{theorem}\label{t3}
Set $q>1$. Assume the payoff function $g$ satisfies (\ref{4}) and that conditions 1., 2., 3. and 4. of Theorem \ref{t2} are satisfied. Then : 
\par If $\lim_{n\to+\infty}I_n(q)=I(q)<+\infty$, the limiting model has an equivalent martingale measure which minimises $f^q$-divergence. Furthermore, 
\begin{equation}\label{t41}
\lim_{n\to+\infty}\mathbb{C}_{T}^n=\mathbb{C}_T
\end{equation}
Otherwise, the limiting model may or may not have an equivalent martingale measure with minimal $f^q$-divergence. 
\par If $\lim_{n\to+\infty}I_n(q)=I(q)+a, a>0$ and $g$ is bounded, 
\begin{equation}\label{t42}
\lim_{n\to+\infty}\mathbb{C}_{T}^n=E_{P^*}(g(S))
\end{equation}
where $P^*$ is a measure equivalent to $P$ under which $S$ is not a martingale. 
\par If $\lim_{n\to+\infty}I_n(q)=+\infty$, and $g$ is bounded,  
\begin{equation}\label{t43}
\lim_{n\to+\infty}\mathbb{C}_{T}^n=E_{P}(g(S))
\end{equation}
where $P$ is the initial measure. 
\end{theorem}

\begin{rem} If we drop the assumption $\nu(\{m\})=\nu(\{M\})=0$, the measure $Q$ associated with the Girsanov parameters $(\beta,Y)$ will still be a martingale measure for the limiting model but may only be absolutely continuous with respect to $P$ as we no longer necessarily have $Y>0$ $\nu$-a.s.  
\end{rem}

\section{Proof of Theorems \ref{t1} and \ref{t1a}}

The proof of Theorem \ref{t1} is given in several steps. We will
assume for simplicity and without loss of generality that
$r_n=r=0$. Such an assumption is equivalent to introducing the
processes $\tilde{X}^n$ and $\tilde{X}$ with
$\tilde{X}^n_t=X^n_t-r_nt$, and $\tilde{X}_t=X_t-rt$, $t\geq 0$ or to
replacing $b_n$ by $b_n-r_n$ and $b$ by $b-r$.

\begin{lemma}\label{10a}
Assume the measures $Q^n$ and $Q$ are equivalent to $P^n$ and $P$ respectively and are martingale measures. We denote by $(\beta^n,Y^n)$ and $(\beta,Y)$ the Girsanov parameters of these measures. Under the assumptions of Theorem \ref{t1}, we have 
\begin{equation}
\label{10aa}
\mathcal L(X^n |Q^n) \rightarrow \mathcal L (X |Q).
\end{equation}
\end{lemma}
\noindent \it Proof \rm We will use Theorem VII.2.9 in \cite{JS} which gives necessary and sufficient conditions for (\ref{10aa}) in terms of characteristics. 
\par For this, we use the fact that according to the Girsanov theorem
(cf. Th. III.3.24 p. 159 in \cite{JS}), the characteristics of $X^n$
under $Q^n$ are equal to 
 \\
\begin{equation}\label{10bb}
\left\{
\begin{array}{l}
b^{Q^n} = b_n  + \beta ^n c_n + \int _{\mathbb{R}^*} h(x)\, (Y^n(x) - 1)\nu _n(dx) ,\\
c^{Q_n} = c_n,\\
\nu ^{Q_n}(dx)= Y^n(x) \nu_n(dx)
\end{array}\right.\end{equation}
 where $\beta^n\in\mathbb{R}$ and $Y^n$ is a positive Borel function which satisfies the condition 
 $$\int_{\mathbb{R}^*}|h(x)(Y^n(x)-1)|\nu_n(dx)<+\infty.$$
 The constant $\beta^n$ is not arbitrary since the process $S^n$ is a martingale with respect to $Q^n$. As is well known, in order for $S^n$ to be a martingale, it is necessary and sufficient that
 (cf.\cite{JS}, p. 556)
 \begin{equation}\label{10c}
b_n  +\beta ^n c_n + \frac{1}{2} c_n + \int _{\mathbb{R}^*}[(e ^x - 1) Y ^n(x) - h(x)] \nu _n(dx) = 0
\end{equation}

Then from (\ref{10bb}) and (\ref{10c}) we have :
\begin{equation}\label{bqn} 
b ^{Q^n} = -\frac{1}{2} c_n - \int _{\mathbb{R}^*}[(e ^x - 1 - h(x)) Y ^n(x) ] \nu _n(dx) 
\end{equation}
In the same way, we obtain the characteristics of the process $X$ with respect to $Q$:
$$\left\{
\begin{array}{l}\label{10b}
b^{Q} = b + \beta  c + \int _{\mathbb{R}^*} h(x)\, (Y(x) - 1) \nu (dx) ,\\
c^{Q} = c,\\
\nu ^{Q}(dx)= Y(x) \nu (dx)
\end{array}\right.$$
According to Theorem VII.2.9 p. 355 in \cite{JS}, we have to check the following conditions: 
\begin{itemize}
\item [j)] $\lim_{n\rightarrow\infty} b^{Q ^n} = b^{Q}$
\item [jj)] $\lim_{n\rightarrow\infty} c^{Q ^n} + \int
  _{\mathbb{R}^*} h ^2(x) \nu ^{Q ^n}(dx) =  c^{Q } + \int _{\mathbb{R}^*} h^2(x) \nu ^{Q }(dx)$
\item  [jjj)] for all continuous bounded functions $f$ such that $\lim_{x\to 0}\frac{f(x)}{x^2}=0$, 
$$\lim_{n\rightarrow\infty}\int _{\mathbb{R}^*} f(x) \nu ^{Q ^n}(dx) = \int _{\mathbb{R}^*} f(x) \nu ^{Q }(dx) $$
\end{itemize}
The conditions jj) and jjj) follow directly from conditions 2. and 3. of Theorem \ref{t1}, since $\nu^{Q^n}=Y^n\nu_n$ and $\nu^Q=Y\nu$. It remains to check j). 
We have from (\ref{bqn}) :
$$b ^{Q ^n} = -\frac{1}{2} [ c_n + \int _{\mathbb{R}^*} h^2(x) Y ^n(x)\nu _n(dx) ] - \int _{\mathbb{R}^*} [ e ^x -1 -h(x) - \frac{h ^2 (x)}{2}]Y^n(x) \nu _n(dx)$$
It follows from conditions 2. and 3. applied to $f(x)=[e^x-1-h(x)-\frac{h^2(x)}{2}]/(e^x-1)$ that both terms on the right-hand side of this equality converge. Hence, the conditions of Theorem VII.2.9 are satisfied and we have (\ref{10aa}). \QED

\begin{lemma}\label{L2}
Assume the processes $X^n$ and $X$ are Levy processes with respect to some equivalent martingale measures $Q^n$ and $Q$ respectively and that 
$$\mathcal{L}(X^n|Q^n)\longrightarrow \mathcal{L}(X|Q).$$
Then for any random variable $\tau$, independent of $X^n$ and $X$ and with an exponential distribution $\mu_q$ of parameter $q>0$, we have 
\begin{equation}
\label{10d}
\lim_{n\rightarrow \infty} E_{Q ^n \times \mu _q} (\sup_{0\leq t \leq \tau} e ^{X_t^n}) =
 E_{Q \times \mu _q} (\sup_{0\leq t \leq \tau}
e ^{X_t}).  
\end{equation}
\end{lemma}
\noindent \it  Proof \rm We  assume that $\tau$ is a random variable
with an exponential distribution of parameter $q>0$ given on a space $(E,\mathcal{E})$. We consider an enlargement of the initial probability space $(\tilde{\Omega}^n,\tilde{\mathcal{F}}^n,\tilde{\mathbb{F}}^n,\tilde{Q}^n)$ with $\tilde{\Omega}^n=\Omega^n\times E$, $\tilde{\Omega} ^n= \Omega ^n \times E$, $\tilde{\mathcal F} ^n= \mathcal F ^n\times \mathcal E$, $\tilde{\mathbb F} ^n = (\tilde{\mathcal F}_t ^n )_{t\geq 0}$ and $\tilde{\mathcal F}_t^n = \mathcal F_t^n\otimes \mathcal E$ и $\tilde{Q} ^n = Q \times \mu _q$. We define in the same way an enlargement $(\tilde{\Omega}, \tilde{\mathcal F}, \tilde{\mathbb F}, \tilde{Q})$ of $(\Omega,\mathcal{F},\mathbb{F},Q)$. 
The processes $X^n$ and $X$ remain Levy processes with the same characteristics under this enlargement.   
 \par According to a result of Rogozin \cite{R}, under the condition $E_{Q^n}[e^{X^n_1}]<e^q$, and on the set $\{z\in\mathbb{C}|Re(z)<1\}$, we have the Wiener-Hopf factorisation 
$$ \frac{q}{q- \ln (E_{Q^n }e ^{zX _1^n})} = E_{\tilde{Q} ^n}[e ^{z \sup_{0\leq t\leq\tau}X_t ^n}] E_{\tilde{Q} ^n}[e ^{z \inf_{0\leq t\leq\tau}X_t ^n}]$$
Using the fact that $Q^n$ is a martingale measure, and hence $E_{Q^n}[e^{X^n_1}]=E_{Q^n}[e^{X^n_0}]=1$, we obtain that this decomposition can be extended to $z=1$ so that
\begin{equation}
\label{11}
E_{\tilde{Q} ^n}[e ^{ \sup_{0\leq t\leq\tau}X_t ^n}]  = \frac{1}{E_{\tilde{Q} ^n}[e ^{\inf_{0\leq t\leq\tau}X_t ^n}]}
\end{equation}
Since the processes $X^n$ and $X$ have no fixed points of discontinuity and since $\tau$ is an independent random variable, we have 
$$\mathcal L ( \inf_{0\leq t\leq \tau} X_t ^n\,|\, \tilde{Q }^n) \rightarrow \mathcal L ( \inf_{0\leq t\leq \tau} X_t \, |\,\tilde{Q })$$
Since $e ^{\inf_{0\leq t\leq \tau} X ^n_t} \leq e ^{X_0^n } \leq 1$, we deduce from the Lebesgue dominated convergence theorem that 
\begin{equation}
\label{12}
\lim_{n\rightarrow \infty} E_{\tilde{Q} ^n}[e ^{\inf_{0\leq t\leq\tau}X_t ^n}] = E_{\tilde{Q} }[e ^{\inf_{0\leq t\leq\tau}X_t }]
\end{equation}
Finally, (\ref{11}) and (\ref{12}) give (\ref{10d}).
\QED

\begin{lemma}\label{L3}
Assume that (\ref{10d}) holds. Then for any $T>0$, the family of random variables $(\sup _{0\leq t\leq T}S_t^n)_{n\geq 1}$ is uniformly integrable.
\end{lemma}
\noindent  \it  Proof \rm We recall that
$$\sup _{0\leq t\leq T}S_t^n = S_0^n e^{\sup _{0\leq t\leq T}X_t^n }$$
For any fixed $q>0$, we have
$$e^{\sup _{0\leq t\leq T}X_t^n } \leq e ^{qT}\int _T^{\infty} q e ^{-q u}e^{\sup _{0\leq t\leq u}X_t^n } d u 
\leq  e ^{qT}\int _0^{\infty} q e ^{-q u}e^{\sup _{0\leq t\leq u}X_t^n } d u $$
so we must have 
\begin{equation}
\label{133}
\sup _{0\leq t\leq T}S_t^n \leq e^{q T}E_{\mu_q}(\sup _{0\leq t\leq \tau }S_t^n)
\end{equation}
Using Fubini's theorem, we deduce from (\ref{10d}) that the family of positive random variables $(E_{\mu_q}(\sup _{0\leq t\leq \tau }S_t^n))_{n\geq 1}$ is uniformly integrable. (\ref{133}) then implies that the family $(\sup_{0\leq t\leq \tau}S^n_t)_{n\geq 1}$ is itself uniformly integrable. \QED

\noindent \it Proof of Theorem \ref{t1} \rm
By Lemma \ref{10a} and condition 1. in Theorem \ref{t1}, we have 
$$\mathcal L (S^n | Q ^n) \rightarrow \mathcal L( S | Q).$$
and since the functional $g$ is continuous on $D[0,T]$, we also have 
$$\mathcal{L}(g(S^n)|Q^n)\rightarrow \mathcal{L}(g(S)|Q).$$
Lemma \ref{L3} then implies convergence (\ref{5}) of option prices for payoff functions which satisfy (\ref{4}). \QED

\noindent \it Proof of Theorem \ref{t1a} \rm  It isn't difficult to show using (\ref{f}) that 
$$\int_{\mathbb{R}^*}(\sqrt{Y(x)}-1)^2\nu(dx)<+\infty$$
and hence that the Hellinger process $h_{T}(P^*,P)$ of order $1/2$ is finite. This means (cf. \cite{JS} Th. IV.2.1 p. 209) that $P^*<<P$. We now need to check that the characteristics of the process $X^n$ with respect to $Q^n$ given in (\ref{10bb}) satisfy j), jj) and jjj). 
We can write
$$\begin{aligned}
b^{Q^n}=&b_n+[c_n+\int_{\mathbb{R}^*}h^2(x)\nu_n(x)]\beta +c_n(\beta^n-\beta)
\\&+\int_{\mathbb{R}^*}h(x)(Y^n(x)-Y(x))\nu_n(dx)+\int_{\mathbb{R}^*}h(x)[Y(x)-1-\beta
h(x)]\nu_n(dx)
\end{aligned}$$
Conditions (\ref{a}),(\ref{e}),(\ref{f}) and i) then give 
$$\lim_{n\to+\infty}b^{Q^n}=b+c\beta+\int_{\mathbb{R}^*}h(x)(Y(x)-1)\nu(dx)$$
In the same way, we show that 
$$\lim_{n\to+\infty}c^{Q^n}+\int_{\mathbb{R}^*}h^2(x)\nu^{Q^n}(dx)=c+\int_{\mathbb{R}^*}h^2(x)Y(x)\nu(dx)$$
and that for every continuous bounded function $f$ which satisfies $\lim_{x\to 0}\frac{f(x)}{x^2}=0$, 
$$\lim_{n\to+\infty}\int_{\mathbb{R}^*}f(x)\nu^{Q^n}(dx)=\int_{\mathbb{R}^*}f(x)Y(x)\nu(dx)$$
Hence, conditions j), jj) and jjj) are satisfied so that by Theorem VII.2.9 in \cite{JS}, we have 
$$\mathcal{L}(X^n|Q^n)\longrightarrow \mathcal{L}(X|P^*).$$
The first part of the theorem then follows from the fact that $\lim_{n\to+\infty}S^n_0=S_0$ and from the continuity and boundedness of $g$ on $D([0,T])$. 
\par It is well known that $P^*$ will be an equivalent martingale measure for $S$ if and only if we have $Y>0$ ($\nu$-a.s) and 
$$b+\frac{c}{2}+\beta c+\int_{\mathbb{R}^*}[(e^x-1)Y(x)-h(x)]\nu(dx)=0$$
Using the Girsanov formula for $b^{P^*}$, this yields 
\begin{equation}\label{bpp}
b^{P^*}=-\frac{c}{2}-\int_{\mathbb{R}^*}(e^x-1-h(x))Y(x)\nu(dx)
\end{equation}
Since $b^{P^*}=\lim_{n\to+\infty}b^{Q^n}$, we can deduce from the
assumptions of our theorem that 
$$b^{P^*}=-\frac{c}{2}-\int_{x<1}(e^x-1-h(x))Y(x)\nu(dx)-\lim_{n\to+\infty}\int_{x\geq 1}(e^x-1)Y^n(x)\nu_n(dx)$$
But (\ref{bpp}) then becomes 
$$\lim_{n\to+\infty}\int_{x\geq 1}(e^x-1)Y^n(x)\nu_n(dx)=\int_{x\geq1}(e^x-1)Y(x)\nu(dx)$$
Hence, the conditions of Theorem \ref{t1} are satisfied and we have convergence of the option prices. \QED

\section{Proof of Theorem \ref{t2}}

We may assume as before without loss of generality that $r_n=r=0$. We denote by $\hat{\psi}^n$ and $\hat{\psi}$ the characteristic exponents of the processes $\hat{X}^n$ and $\hat{X}$ which are used to define minimal entropy martingale measures. We also consider the sets
$$D_n = \{ u\in \mathbb R : \int_{x>1} e ^{u(e ^x-1)} \nu _n (dx) < \infty \}$$
and 
$$D = \{ u\in \mathbb R : \int_{x>1} e ^{u(e ^x-1)} \nu (dx) < \infty \}$$
on which the Laplace transforms of $\hat{X}^n_1$ and $\hat{X}_1$ are defined. Then according to \cite{MI}, \cite{MIF}, for each $u\in D_n$,
\begin{equation}
\label{13}
\hat{\psi}_n(u) = (b_n +\frac{c_n}{2})u + \frac{c_n}{2} u ^2 + \int_{\mathbb{R}^*} \left (  e ^{u(e ^x-1)} - 1 - u h(x)\right) \nu_n(dx)
\end{equation}
and the same representation holds for $\hat{\psi}$ : if $u\in D$,
$$\hat{\psi}(u) = (b +\frac{c}{2})u + \frac{c}{2} u ^2 + \int_{\mathbb{R}^*} \left (  e ^{u(e ^x-1)} - 1 - u h(x)\right) \nu(dx)$$
First we recall the following important fact :
\begin{lemma}\label{L30}(cf \cite{ES},\cite{MIF})
We assume that there exists an equivalent martingale measure $Q^n$ with minimal entropy for $X^n$. The Girsanov parameters $(\beta^n,Y^n)$ are deterministic functions and the process $X^n$ remains a Levy process under $Q^n$. 
\end{lemma}
\noindent We also recall that we have assumed that $X$ is not a monotone Levy
process and that a Levy process is monotone if and only if we are in one of the following situations :
\par (i)$c=0,\int_{\mathbb{R}^*}h(x)\nu(dx)<+\infty,\nu(\mathbb{R}^{-,*})=0$ and $b-\int_{\mathbb{R}^*}h(x)\nu(dx)\geq 0$
\par (ii) $c=0, \int_{\mathbb{R}^*}h(x)\nu(dx)<+\infty, \nu(\mathbb{R}^{+,*})=0$ and $b-\int_{\mathbb{R}^*}h(x)\nu(dx)\leq 0$. 
\par \noindent It is then easy to show that 
\begin{lemma}\label{L31}
$X$ is a monotone Levy process if and only if
$$\lim_{u\to-\infty}\hat{\psi}'(u)\geq 0 \text{ or }\lim_{u\to+\infty}\hat{\psi}'(u)\leq 0$$
\end{lemma}

We now consider the convergence of the functions $\hat{\psi}^n$. On the set $\cup_{n=1}^{\infty} \cap_{k\geq n} D_k$, the characteristic exponents are defined for large values of $n$, and we can consider their limits.
The set $U$ defined in (\ref{u}) then satisfies 
$$U \subseteq \cup_{n=1}^{\infty} \cap_{k\geq n} D_k.$$
\begin{lemma}\label{L32}
Under the assumptions of Theorem \ref{t2}, the functions $\hat{\psi}^n$ and $\hat{\psi}$ belong to $C^{\infty}(\stackrel{\circ}{D^n})$ and $C^{\infty}(\stackrel{\circ}{D})$ respectively. Moreover, 
\begin{equation}\label{14}
\hat{\psi}_n'(u) \rightarrow \hat{\psi}'(u) 
\end{equation}
on any compact set $K\subset \stackrel{\circ}{U}$ as $n\rightarrow +\infty$.
\end{lemma}
\noindent \it Proof \rm We prove that $\hat{\psi} \in
C^{\infty}(\stackrel{\circ}{D})$. The fact that $\hat{\psi}^n \in
C^{\infty}(\stackrel{\circ}{D})$ can obviously be shown in the same
way. Note that $\stackrel{\circ}{D}$ is an open interval. Since $h(x)=x$ in a neighbourhood of $0$ and using the mean-value theorem, there exists a constant $c>0$ such that for $|x|\leq 1$
$$|\,e ^{u (e ^x - 1)} - 1 - uh(x)\,| \leq c(u ^2 + u)x ^2$$
The uniform convergence of the integrals on $\{|x|<1\}$ then follows from (\ref{1}) (where $\nu_n$ is replaced by $\nu$). If $|x|\geq 1$, as $\nu([1,+\infty[)<+\infty$ and $h(x)=0$, it is sufficient to consider the integral 
$$\int_{|x|>1}e ^{u (e ^x - 1)} \nu (dx) =  \int_{x<-1}e ^{-u (1-e ^x )} \nu (dx)+\int_{x>1}e ^{u (e ^x - 1)} \nu (dx) $$
Both integrals on the right-hand side converge uniformly in $u$, as
the integrating functions are uniformly bounded from above
respectively by $e^{|u|}$ and $e^{(u+\delta)(e^x - 1)}$  where $\delta
>0$ is such that $u+\delta\in \stackrel{\circ}{D}$. Similar considerations can be applied to the derivatives $\hat{\psi}^{(k)},\,\, k\geq 1$.
\par In order to prove the uniform convergence given in (\ref{14}), we note that for $u\in K$, 
\begin{equation}\label{15}
\hat{\psi}'(u) = b+(\frac{1}{2}+u) c + \int_{\mathbb{R}^*} \left ( (e^x - 1) e ^{u(e ^x-1)} - h(x)\right) \nu(dx)
\end{equation}
Similar formulae hold for $\hat{\psi}'_n(u)$. 
\par Since $\hat{\psi}'_n$ and $\hat{\psi}'$ are continuous increasing functions on $K$, uniform convergence of the sequence is equivalent to point-wise convergence. In view of conditions 2.,3. and (\ref{15}) it is enough to show  that for $u\in K$
$$\lim_ {n\rightarrow \infty} \int _{\mathbb{R}^*} f (u,x) \nu _n(dx) = \int _{\mathbb{R}^*} f(u,x) \nu (dx) $$
where $f(u,x)= (e^x - 1) e ^{u(e ^x-1)}- h(x) - h ^2(x)(u+
\frac{1}{2})$. This follows from the fact that $\lim_{x\to
  0}\frac{f(u,x)}{x^2}=0$ and that since $u\in K\subseteq
\stackrel{\circ}{D}$, the family of functions is  uniformly integrable.\QED

\begin{lemma}\label{L33}
Assume the conditions of Theorem \ref{t2} are satisfied. If $\lim_{u\to\alpha^-}\hat{\psi}'(u)>0$, equation (\ref{8}) has a solution $\theta<\alpha$;  
if $\lim_{u\to\alpha^-}\hat{\psi}'(u)=0$, then $\theta=\alpha$ with
$\alpha\neq \infty$; in both cases (in the last case for a subsequence)  
\begin{equation}\label{17}
\lim_{n\rightarrow \infty} \theta _n = \theta.
\end{equation}
If $\lim_{u\to\alpha^-}\hat{\psi}'(u)<0$, (\ref{8}) may or may not have a solution, but
there exists a subsequence $n'$ such that
\begin{equation}\label{18}
\lim_{n'\to+\infty} \theta _{n'}  = \alpha.
\end{equation}
 \end{lemma}
\noindent \it Proof \rm
Assume $\lim_{u\to\alpha^-}\hat{\psi}'(u)>0$. Since $\hat{\psi}'$ is continuous on $\stackrel{\circ}{U}\,\subset\, \stackrel{\circ}{D}$, and it follows from Lemma \ref{L31} that $\lim_{u\to-\infty}\hat{\psi}'(u)<0$, there exists a solution $\theta<\alpha$ to (\ref{8}). 
As $\hat{\psi}'_n$ is an increasing function, we also have that
\begin{equation}
\label{18a}
\lim_{n\rightarrow\infty}\lim_{u\rightarrow \alpha -} \hat{\psi}'_n(u) \geq \lim_{u\rightarrow \alpha -} \hat{\psi}'(u) .
\end{equation}
so that for $n\geq n_0$, $\lim_{u\to\alpha^-}\hat{\psi}_n'(u)>0$ and $\theta_n<\alpha$. Therefore (\ref{17}) follows from the uniform convergence obtained in Lemma \ref{L32}. It should be noted that two facts have played an important part here : first that the solution to equation (\ref{8}) is unique and secondly that $X$ is not a monotone Levy process. 
\par We now assume that $\lim_{u\to\alpha^-}\hat{\psi}'(u)\leq 0$. Since $X$ is not monotone Levy process, $\alpha\neq \infty$. For $\alpha< \infty$, using Fatou's lemma, 
$$\int_{x>1}(e ^x -1) e ^{\alpha(e ^x-1)} \nu(dx)  \leq \underline {\lim}_{u\rightarrow \alpha -}\int_{x>1} (e ^x -1) e ^{u(e ^x-1)} \nu (dx) < \infty$$
Therefore, $\hat{\psi}'(\alpha)$ is well defined and from the Lebesgue
monotone convergence theorem we have
$\hat{\psi}'(\alpha)=\lim_{u\to\alpha^-}\hat{\psi}'(u)$. In
particular, if this limit is $0$, the equation $\hat{\psi}'(u)=0$ has
a solution $\theta=\alpha$. In both cases, for $u>\alpha$, it follows
from the definition of $U$ given in (\ref{u}) that there is a
subsequence $n'$ for which
$\lim_{n'\to+\infty}\hat{\psi}'_{n'}(u)=+\infty$. Hence, for $n$ large
enough $\theta_{n'}<u$ and $\lim_{n'\to+\infty}\theta_{n'}\leq \alpha$. Furthermore, for every $u<\alpha$, $\hat{\psi}'(u)<0$ and hence by Lemma \ref{L32}, $\hat{\psi}_n'(u)<0$ for $n$ big enough so that $\theta_n>u$. Thus $\lim_{n'\to+\infty}\theta_{n'}=\alpha$. \QED
 
\noindent \it Proof of Theorem \ref{t2}  \rm If $\lim_{u\to\alpha^-}\hat{\psi}'(u)\geq 0$, we can apply both Theorems \ref{t1} and \ref{t1a}. We show for example that the conditions of Theorem \ref{t1a} hold. Since $Y^n(x)=e^{\theta_n (e^x-1)}$, $Y(x)=e^{\theta (e^x-1)}$, and (\ref{17}) holds, we have 
\begin{equation}\label{20a}
Y^n(x) = 1 + \theta _n x + o(x) \text{ and }Y(x) = 1 + \theta  x + o(x)
\end{equation}
where $o(x)$ does not depend on $n$. (\ref{17}) gives condition i) of
Theorem \ref{t1a}. 
\par Since $\hat{\psi^n}'(\theta_n)=0$, we have
\begin{equation} \label{defthetan}
\int_{x\geq \epsilon}(e^x-1)e^{\theta_n
  (e^x-1)}\nu_n(dx)=-[b_n+c_n(\frac{1}{2}+\theta_n)+\int_{x<\epsilon}[(e^x-1)e^{\theta_n(e^x-1)}-h(x)]\nu_n(dx)]
\end{equation}
It follows from 2., 3. and 4. and (\ref{17})  that the
expressions on the right-hand side are uniformly bounded.
This implies
\begin{equation}\label{ui}
\sup_n \int_{x\geq\epsilon}(e^x-1)e^{\theta_n (e^x-1)}<+\infty
\end{equation}
This, and the fact that on the set $\{x\leq-\epsilon\}$ the functions $(e^{\theta_n (e^x-1)})$ are
bounded by the constant $e^{\theta +1}$, implies that the family
$(e^{\theta_n (e^x-1)})_{n\geq 1}$ is uniformly integrable on the set
$\{|x|\geq\epsilon \}$, so that 
\begin{equation}\label{cdii}
\lim_{n\to+\infty}\int_{|x|\geq \epsilon}e^{\theta_n (e^x-1)}\nu_n(dx)=\int_{|x|\geq \epsilon}e^{\theta (e^x-1)}\nu(dx)
\end{equation}
Furthermore, it follows from (\ref{defthetan}) and $\hat{\psi}'(\theta)=0$
that 
\begin{equation}\label{cdiii}
\lim_{n\to+\infty}\int_{x\geq 1}(e^x-1)e^{\theta_n
  (e^x-1)}\nu_n(dx)=\int_{x\geq 1}(e^x-1)e^{\theta (e^x-1)}\nu(dx)
\end{equation}
and so $P^*$ is a martingale measure. 
\par If now $\lim_{u\to\alpha}\hat{\psi}'(u)<0$, we have
(\ref{18}), and (\ref{20a}) and (\ref{cdii}) remain true with $\alpha$ instead of
$\theta$. However, (\ref{cdiii}) no longer holds and $P^*$ is not a
martingale measure. 
\par In all cases, the conditions of Theorem \ref{t1a} are satisfied
and Theorem \ref{t2} is proved. \QED

\begin{exple}\rm Let for every $n\geq 1$ $V^n$ be an NIG$(n,0,n,0)$ process and $Z^n$ be an NIG$(\frac{1}{4},0,\frac{1}{n},0)$ process, independent of $V^n$.  We recall that the Levy measure of a NIG process satisfies (\ref{levnig}). We consider the sequence of Levy processes $(X^n)_{n\geq 1}$ defined by 
$$X^n_t=bt+V^n_t+Z^n_t$$ 
so that
 $$\hat{\psi}_n(u)=ub+\frac{1}{\pi}\int_{\mathbb{R}^*}(e^{u(e^x-1)}-1-uh(x))\frac{n^2 K_1(n|x|)+\frac{1}{4n}K_1(\frac{|x|}{4})}{|x|}dx$$
 where $K_1$ is the modified Bessel function of the third kind of parameter 1. Its behaviour at $0$ and $+\infty$ is given by (cf. formulae 9.6.9 and 9.7.2 in \cite{AS})
 $$K_{1}(z)\sim \frac{1}{z} \text{ when }z\longrightarrow 0 \text{ and }K_1(z)\sim \sqrt{\frac{\pi}{2z}}e^{-z} \text{ when }z\longrightarrow +\infty$$
 For every $n$, $\hat{\psi}_n$ is defined on $]-\infty,0]$ and differentiable on $]-\infty,0[$ so that $\alpha=0$. We can check that 
 $$\lim_{u\to-\infty}\hat{\psi}_n'(u)=-\infty\text{ and }\lim_{u\to 0}\hat{\psi}_n'(u)=+\infty$$
 so that there exists for each $S^n$ an equivalent martingale measure of minimal entropy. As $n$ goes to $+\infty$, 
 $$\mathcal{L}(X^n|P^n)\longrightarrow (bt+W_t)_{t\geq 0}$$
 where $W$ is a standard Brownian motion. Thus, we have $\lim_{u\to 0}\hat{\psi}'(u)=b+\frac{1}{2}$. Applying Theorem 3, if $b+\frac{1}{2}\geq 0$, we have convergence to the price in the Black-Scholes formula, i.e.
$$\lim_{n\to+\infty}\mathbb{C}_T^n=E_{Q}[g(S)]$$ 
where $Q$ is the unique martingale measure for the limiting  model.
If now $b+\frac{1}{2}<0$ and $g$ is bounded, 
$$\lim_{n\to+\infty}\mathbb{C}_T^n=E_{P}[g(S)]$$ where $P$ is the initial measure. The limiting model  
does however have a unique equivalent martingale measure $Q\neq P$. In particular, if we consider a european put option with maturity $T$ and exercise price $K$, we have 
$$\lim_{n\to+\infty}E_{Q^n}(K-S_T^n)^+=E_P(K-S_T)^+>E_Q(K-S_T)^+$$ 
\end{exple}

\begin{exple} \rm We now consider a sequence of NIG$(\alpha_n,-\alpha_n,1,-1)$ processes with $\alpha_n=\frac{1}{2}-\frac{1}{4n}$. We then have 
$$\hat{\psi}_n(u)=-u+\frac{\alpha_n}{\pi}\int_{\mathbb{R}^*}(e^{u(e^x-1)}-1-uh(x))e^{-\alpha_n x}\frac{K_1(\alpha_n |x|)}{|x|}dx$$
The function $\hat{\psi}_n$ is defined on $]-\infty,0]$ and
differentiable on $]-\infty,0[$ so that $\alpha=0$ as before. We can
check that $\lim_{u\to-\infty}\hat{\psi}_n'(u)=-\infty$ and that if
$\alpha_n<\frac{1}{2}$, $\lim_{u\to
  0}\hat{\psi}_n'(u)=+\infty$. Therefore, the equation
$\hat{\psi^n}'(u)=0$ has a solution and there is an equivalent
martingale measure $Q^n$ with minimal entropy. We can show furthermore that 
$$\lim_{n\to+\infty}b^n=-1, \lim_{n\to+\infty}\int_{\mathbb{R}^*}h^2(x)\nu_n(dx)=\frac{1}{2\pi}\int_{\mathbb{R}^*}h^2(x)e^{-\frac{x}{2}}\frac{K_1(\frac{|x|}{2})}{|x|}dx$$
and that for every continuous bounded function $f$ which satisfies $\lim_{x\to 0}\frac{f(x)}{x^2}=0$, 
$$\lim_{n\to+\infty}\int_{\mathbb{R}^*}f(x)\nu_n(dx)=\frac{1}{2\pi}\int_{\mathbb{R}^*}f(x)e^{-\frac{x}{2}}\frac{K_1(\frac{|x|}{2})}{|x|}dx$$
Conditions 2., 3. and 4. of Theorem 3 are therefore satisfied for a
limiting process $X$ which is a NIG$(\frac{1}{2},-\frac{1}{2},1,-1)$
process. The function $\hat{\psi}'$ is defined on $]-\infty,0]$ and we have $\hat{\psi}'(0)=-1$. The equation  $\hat{\psi}'(u)=0$ does not have a solution and the limiting model does not have an equivalent martingale measure of minimal entropy. However, using theorem 3, we have for instance for a european put option 
$$\lim_{n\to+\infty}E_{Q^n}(K-S_T^n)^+=E_{P}(K-S_T)^+$$
\end{exple} 

\section{Proof of Theorem \ref{t3}} 

Here again, we can assume for simplicity and without loss of
generality that $r_n=r=0$. 
We set $q>1$ and we assume as before that $X$ is not a monotone Levy process. Assuming $I(q)<+\infty$, we  introduce the function 
\begin{equation}\label{defF}
F(u)= b+\frac{c}{2}+cu+\int_{\mathbb{R}^*}[(e^x-1)Y_u(x)-h(x)]\nu(dx)
\end{equation}
where 
\begin{equation}\label{defyu}
Y_u(x)=\begin{cases}[1+(q-1)u(e^x-1)]^{\frac{1}{q-1}}&\text{ if }(q-1)u(e^x-1)\geq -1,
\\0 &\text{ otherwise, }
\end{cases}
\end{equation}
and we can show as for Lemma \ref{L31} :
\begin{lemma}\label{mon}
$X$ is monotone if and only if 
$$\lim_{u\to-\infty}F(u)\geq 0 \text{ or } \lim_{u\to+\infty}F(u)\leq 0.$$
\end{lemma}
We first consider the set of martingale measures which are absolutely
continuous with respect to $P$ and we show that the Girsanov
parameters of the minimal $f^q$-martingale measure are the solutions
to a minimisation problem. 

\begin{lemma}\label{min} The Girsanov parameters $(\beta,Y_{\beta})$ of the absolutely continuous martingale measure $Q$ which minimises the Hellinger integral of order $q$ are deterministic functions. If $c\neq 0$, $Y_{\beta}$ is the unique solution to the problem of minimising the function 
$$k(Y)=\frac{1}{2}\frac{q(q-1)}{c}[(b+\frac{c}{2}+\int_{\mathbb{R}^*}(e^x-1)Y(x)-h(x)\nu(dx))]^2+\int_{\mathbb{R}^*}(Y^q(x)-q(Y(x)-1)-1)\nu(dx)$$
on the Banach space of piecewise continuous non-negative functions with the locally uniform norm. 
\par If $c=0$, the parameter $Y_{\beta}$ is the unique solution for the minimisation of $k_0(Y)$ :
$$k_0(Y)=\int_{\mathbb{R}^*}(Y^q(x)-q(Y(x)-1)-1)\nu(dx)$$
under the constraint $b+\int_{\mathbb{R}^*}(e^x-1)Y(x)-h(x)\nu(dx)=0$. In
both cases, $Y_{\beta}$ is given by (\ref{defyu}) with $u=\beta$,
where $\beta$ is the unique solution to the equation $F(u)=0$. 
\end{lemma}

\it Proof \rm It follows from Lemma \ref{L31} that as $X$ is not
monotone, the equation $F(u)=0$ has a solution which is furthermore
unique since $F'(u)>0$ on $\mathbb{R}$. Following the proof of Theorem 2.9 in \cite{JKM}, we see that the martingale measure with Girsanov parameters $(\beta,Y_{\beta})$  minimises the Hellinger integral of order $q$. \QED
\vskip 0.1cm
\par For every $n\geq 1$, we define an analogue to function $F$ :
$$F^n(u)= b_n+\frac{c_n}{2}+c_nu+\int_{\mathbb{R}^*}[(e^x-1)Y_u(x)-h(x)]\nu_n(dx)$$
It follows from Lemma \ref{min} that the Girsanov parameters $(\beta^n,Y^n)$ of the minimal measure $Q^n$ satisfy $F^n(\beta_n)=0$ and $Y^n=Y_{\beta^n}$. 

\begin{lemma}\label{cvpsin}
If 
\begin{equation}\label{alpha}
\lim_{n\to+\infty}I_n(q)=I(q)<+\infty,
\end{equation} then $F^n$ converges uniformly on compact sets to $F$. If 
\begin{equation}\label{beta}
\lim_{n\to+\infty}I_n(q)=I(q)+a, a>0,
\end{equation}
 then $F^n$ converges uniformly on compact sets to $\tilde{F}>F$. If 
\begin{equation}\label{gamma}
\lim_{n\to+\infty}I_n(q)=+\infty,
\end{equation} then for every $u>0$, $\lim_{n\to+\infty}F^n(u)=+\infty$.
\end{lemma}
\vskip 0.2cm 
\noindent \it Proof \rm For every $u\in\mathbb{R}$, conditions 2., 3. and 4. imply 
\begin{equation}\label{ast1}
\lim_{n\to+\infty}b_n+\frac{c_n}{2}+c_n u+\int_{x<1}[(e^x-1)Y_u(x)-h(x)]\nu_n(dx)=b+\frac{c}{2}+cu+\int_{x<1}[(e^x-1)Y_u(x)-h(x)]\nu(dx)
\end{equation}
where $Y_u$ is defined by (\ref{defyu}).  
\par Now there exists a positive constant $C$ such that for $x\geq 1$,  
\begin{equation}\label{ko}
(e^x-1)Y_u(x)\leq C e^{\frac{qx}{q-1}}
\end{equation}
If (\ref{alpha}) holds, the function $(e^x-1)Y_u(x)$ is
uniformly integrable with respect to $\nu_n$ and
$$\lim_{n\to+\infty}\int_{x\geq 1}(e^x-1)Y_u(x)\nu_n(dx)=\int_{x\geq 1}(e^x-1)Y_u(x)\nu(dx)$$
Hence $\lim_{n\to+\infty}F^n(u)=F(u)$. 
\par If however we have (\ref{beta}), then for $u\neq 0$,  
$$\lim_{n\to+\infty}\int_{x\geq 1}(e^x-1)Y_u(x)\nu_n(dx)>\int_{x\geq 1}(e^x-1)Y_u(x)\nu(dx)$$
so that
$\lim_{n\to+\infty}F^n(u)=\tilde{F}(u)>F(u)$.
Furthermore, in both these cases, all the functions we consider are continuous and strictly increasing, so the convergence is uniform on compact sets. 
\par Finally, if we have (\ref{gamma}), then $\nu_n(\{x\geq 1\})>0$ for $n$ large enough. For every $u>0$, $Y_u>0$ on $\{x\geq 1\}$ and we have
$$\lim_{n\to+\infty}\int_{x\geq 1}(e^x-1)Y_u(x)\nu_n(dx)=+\infty$$
and so $\lim_{n\to+\infty}F^n(u)=+\infty$. \QED

\begin{lemma}\label{cvbetan}
If $\lim_{n\to+\infty}I_n(q)<+\infty$, then $\lim _{n\rightarrow \infty}\beta_n=
\beta$ where $\beta$ is the unique solution of the
equation $F(u)=0$ when (\ref{alpha}) holds and of the equation
$\tilde{F}(u)=0$ when (\ref{beta}) holds.  
If $\lim_{n\to+\infty}I_n(q)=+\infty$, then $\lim_{n\to+\infty}\beta_n=0$.
\end{lemma}

\noindent \it Proof \rm We first consider the case (\ref{alpha}). We
assume that $\underline{\beta}=\underline{\lim}_{n\to+\infty}
\beta_n<\overline{\beta}=\underline{\lim}_{n\to+\infty}
\beta_n$. First of all, as $X$ is not monotone, it follows from Lemma
\ref{mon} that $\underline{\beta}\neq -\infty$ and
$\overline{\beta}\neq +\infty$ so that
$[\underline{\beta},\overline{\beta}]$ is a compact set of
$\mathbb{R}$. If we assume that $\underline{\beta}\neq
\overline{\beta}$, then we can find two subsequences
$\beta^{n'}\rightarrow \underline{\beta}$ and $\beta^{n''}\rightarrow
\overline{\beta}$. Due to uniform convergence, this leads to
$F(\underline{\beta})=0$ and $F(\overline{\beta})=0$ which is
impossible since the solution to $F(u)=0$ is unique. 
\par The case (\ref{beta}) can be treated in exactly the same way by replacing $F$ by the function $\tilde{F}$ introduced in Lemma \ref{cvpsin}. 
\par If we now have (\ref{gamma}), then according to Lemma
\ref{cvpsin}, for every $u>0$, $\lim_{n\to+\infty}F^n(u)=+\infty$. 
Hence for all $u>0$, and $n$ large enough, $\beta_n<u$ so that $\overline{\lim}_{n\to+\infty}\beta_n\leq 0$. 
On the other hand,  (\ref{gamma}) also means that
$\lim_{n\to+\infty}\sup supp(\nu_n)=+\infty$. As the measure $Q^n$ is equivalent to $P^n$, we
must have $Y^n>0\,$ $(\nu_n$-a.s.) and so for $n$ large
enough, $1+\beta_n(q-1)(e^x-1)>0$ for
big values of $x$. Hence $\underline{\lim}_{n\to+\infty}\beta_n\geq 0$, so that finally $\lim_{n\to+\infty}\beta_n=0$. \QED

\vskip 0.2cm\noindent\it Proof of Theorem \ref{t3} \rm We verify the
conditions of Theorem \ref{t1a}. According to the previous lemma, in
each case, $\lim_{n\to+\infty}\beta_n=\beta\in\mathbb{R}$. Hence we
can define $Y_{\beta}$ by (\ref{defyu}). Since
$supp(\nu_n)\subseteq \{x: Y^n>0\}$, condition 4. and the assumption
that $\nu(\{m\})=\nu(\{M\})=0$ implies that $supp(\nu)\subseteq \{x :
Y>0\}$. 
\par For every $n\geq 1$, 
$$Y^n(x)=1+\beta_n x+o(x)\text{ and }Y_{\beta}(x)=1+\beta x+o(x)$$
It follows from Lemma \ref{cvbetan} that $\lim_{n\to+\infty}\beta_n=\beta$ which implies in turn that $o(x)$ can be taken to be independent of $n$. 
\par We
see that in the case when $\lim_{n\to+\infty}I_n(q)<+\infty$, the sequence of
functions $(Y^n)_{n\geq 1}$ is uniformly integrable on
$|x|>\epsilon$. In fact, the functions $Y^n$ are bounded by a
constant for $x<-\epsilon$ and for $x>\epsilon$, (\ref{ko}) holds with
$C=C(\epsilon)$.  
\par Furthermore, if (\ref{alpha}) holds, we also have 
$$\lim_{n\to+\infty}\int_{x\geq 1}e^xY^n(x)\nu_n(dx)=\int_{x\geq 1}e^xY_{\beta}(x)\nu(dx)$$
so that $P^*=Q$ is an equivalent martingale measure and for $g$ which
satisfies (\ref{5}) we have convergence (\ref{t41}) for option prices.  
\par Finally, if we have (\ref{gamma}), $\beta=0$ and
$Y(x)=1$. It follows from the fact that $F^n(\beta_n)=0$ and
conditions 2., 3. and 4. that 
$$\lim_{n\to+\infty}\int_{x\geq \epsilon}(e^x-1)Y^n(x)
\nu_n(dx)=-b-\frac{c}{2}-\int_{x<\epsilon}(e^x-1)-h(x)\nu(dx)$$
It is then easy to see that the family $(Y^nI_{|x|>\epsilon})$ is
uniformly integrable and that condition (ii) of theorem \ref{t1a}
holds. \QED

\begin{rem}
We can obtain results similar to Theorem \ref{t3} when $q<1$. If for
example $0<q<1$, we necessarily have (\ref{alpha}). However,
one has to be a little careful about what happens around the
boundaries of the support of $\nu$, as the sequence $Y^n$ may no
longer be uniformly integrable on $\{|x|>\epsilon\}$ and condition
ii) of Theorem \ref{t1a} may no longer be satisfied. More precisely, the assumptions of Theorem \ref{t3}
no longer ensure the existence of a solution to $F(u)=0$ in
$\mathcal{D}=\{u:1+(q-1)u(e^x-1)>0\ \text{ }\nu-a.s.\}$. If we do have such a solution, the previous
result holds. But if for example $F(u)>0$ on $\mathcal{D}$, one can show that 
$$\lim_{n\to+\infty}\mathbb{C}_{T}^n=E_{P^*}[g(\check{S})]$$
where $\check{S}=exp(\check{X})$ and $\check{X}$ is a Levy process with characteristics $(b,c,\nu+\nolinebreak\frac{F(\inf\mathcal{D})}{1-e^{m}}I_{\{m\}})$. 
\par In the same way, if $F(u)<0$ on $\mathcal{D}$, we have 
$$\lim_{n\to+\infty}\mathbb{C}_T^n=E_{P^*}[(\breve{S})]$$
where $\breve{S}=exp(\breve{X})$ and $\breve{X}$ is a Levy process with characteristics $(b,c,\nu+\frac{|F(\sup\mathcal{D})|}{e^M-1}I_{\{M\}})$.
\end{rem}

\begin{exple} \rm We assume that the processes $X^n$ are
  CGMY$(1,\alpha_n,A_n,B_n)$ processes with drift $b_n$ and parameters
  $ A_n, B_n>0$, $0\leq \alpha_n<1$. We recall that the Levy measure
  of a CGMY process is given by  (\ref{levcgmy}). We consider the case of variance minimising measures (q=2).  We then have  that
\begin{equation}\label{yu2}
Y_u(x)=1+u(e^x-1)
\end{equation}
and 
\begin{equation}\label{ast2}
F^n(u)=b^n+\int_{\mathbb{R}^*}(u(e^x-1)^2+e^x-1-h(x))\nu_n(dx)
\end{equation} 
\par We assume that for every $n\geq 1$, $B^n>2$ so that $F^n$ is well defined on $[0,1]$ and we assume that the equation $F^n(u)=0$ has a solution $\beta_n\in[0,1]$. We also assume that 
$$\lim_{n\to+\infty}b_n=b, \lim_{n\to+\infty}A_n=A, \lim_{n\to+\infty}B_n=B \text{ and }\lim_{n\to+\infty}\alpha_n=\alpha$$
It is easy to see that conditions 2., 3. and 4. are satisfied for a limiting process which is a CGMY$(1,\alpha,A,B)$ process with drift $b$. 
\par If $B>2$, $\lim_{n\to+\infty}I_n(q)=I(q)<+\infty$, so applying Theorem \ref{t3}, 
$$\lim_{n\to+\infty}\mathbb{C}_{T}^n=\mathbb{C}_{T}$$
If however $B=2$, $\lim_{n\to+\infty}I_n(q)=+\infty$. Furthermore, $I(q)=+\infty$, so the limiting model has no variance minimising equivalent martingale measure. But from Theorem \ref{t3} we see that for every bounded payoff function 
$$\lim_{n\to+\infty}\mathbb{C}_{T}^n=E_{P}(g(S))$$ 
\end{exple}

\begin{exple} \rm We assume that $X^n_t=-t+W_t+Z^n_t$ where $W$ is a
  standard Brownian motion and $Z^n$ is a purely discontinuous Levy
  process independent of $W$ for which
  $\nu_n(dx)=\frac{1}{n}e^{-(2+\frac{1}{n})x}I_{[n,2n]}(x)dx$. (\ref{ast2})
  and simple calculations then give
$$F^n(u)=-\frac{1}{2}+\frac{\gamma_n}{n+1}-\frac{\delta_n}{2n+1}+u(\alpha-\frac{2\gamma_n}{n+1}+\frac{\delta_n}{2n+1})$$
where 
$$\alpha=e^{-1}-e^{-2}+1\text{ ,
}\gamma_n=e^{-(n+1)}-e^{-2(n+1)}\text{ , }\delta_n=e^{-(2n+1)}-e^{-2(n+1)}$$
It is easy to see that for every $n$, $F^n(u)=0$ has a solution $\beta_n$ and that the conditions of Theorem \ref{t3} are satisfied for a limiting process $X_t=-t+W_t$. 
\par We have $\lim_{n\to+\infty}I_n(q)=e^{-1}-e^{-2}$ whereas $I(q)=0$. Applying Theorem \ref{t3}, for any continuous bounded payoff function, 
$$\lim_{n\to+\infty}\mathbb{C}_{T}^n=E_{P^*}(g(S))$$
although the limiting model has a unique equivalent martingale measure
$P^*\neq Q$. In particular, for european put options, we obtain an
inequality opposite to that of Example 1 : 
$$\lim_{n\to+\infty}\mathbb{C}_{T}^n<E_{Q}(K-S_T)^+$$
\end{exple}

 \vskip 0.2cm
\noindent \bf{Acknowledgments} \rm Both authors would like to
thank V. Rivero and
M.E Caballero for valuable discussions and suggestions during their stay
in Guanajuato and Mexico as part of the ECOS project M07M01.

\end{document}